%ramacubic.tex: Automatic Solving of Cubic Diophantine Equations Inspired by Ramanujan
%%a Plain TeX file by Shalosh B. Ekhad and Doron Zeilberger (x pages)

%begin macros

\baselineskip=14pt
\parskip=10pt

\font\eightrm=cmr8 

\magnification=\magstephalf

\def\1{{\overline{1}}}
\def\2{{\overline{2}}}
\parindent=0pt
\overfullrule=0in

\def\frac#1#2{{#1 \over #2}}
%\headline={\rm  \ifodd\pageno  \RightHead  \else  \LeftHead  \fi}
%\def\RightHead{\centerline{
%Title
%}}
%\def\LeftHead{ \centerline{Doron Zeilberger}}
%end macros
\centerline
{\bf  
Automatic Solving of Cubic Diophantine Equations Inspired by Ramanujan
}
\bigskip
\centerline
{\it Shalosh B. EKHAD and Doron ZEILBERGER}
\bigskip

{\bf Abstract.} In Ramanujan's Lost Notebook, he stated infinitely many `almost counterexamples' to Fermat's Last Theorem 
for $n=3$, by solving $X^3+Y^3-Z^3=1$. As often the case with Ramanujan, he gave no indication how he discovered it.
Using ingenious {\it exegesis} of Ramanujan's mind by Michael Hirschhorn, combined with 
a much earlier `reading of Ramanujan's mind' by Eri Jabotinsky, we automate this process, and develop
a symbolic-computational algorithm, based on the $C$-finite ansatz, to solve much more general equations,
namely cubic equations of the form $a\,X^3 \,+ \, a\,Y^3 \,+ \, b\,Z^3 \,= \,c$.

{\bf Preface: How this paper was born}

{\bf 2020:} On Sunday, May 24, 2020 we were fortunate to attend a surprise  retirement party for our dear friend Dennis Stanton, and
one of the eight speeches was by the eminent combintorialist and number theorist, George Andrews, who presented
Dennis a retirement present,  the book [Al] {\it ``The  Boy Who Dreamed of Infinity''} authored by his daughter,
Amy Alznauer, meant for ``the children in Dennis' life". As soon as the party was over, we ordered a copy
and after a few days got it. We {\it loved} this book. It turned out that this book can be enjoyed
by {\bf all} ages from $0$ to $\infty$. One of the delightful illustrations, contained the following
{\bf amazing} identity, taken from {\it Ramanujan's Lost Notebook} [R] (p. 341).

{\bf Theorem (Ramanujan)}: Define sequences of integers $a_n,b_n,c_n$, $n \geq 0$ , in terms of the following
generating functions
$$
\sum_{n \geq 0} a_n x^n \, = \, \frac{1+53x+9x^2}{1-82x-82x^2+x^3} \quad ,
$$
$$
\sum_{n \geq 0} b_n x^n \, = \, \frac{2-26x-12x^2}{1-82x-82x^2+x^3} \quad ,
$$
$$
\sum_{n \geq 0} c_n x^n \, = \, \frac{2+8x-10x^2}{1-82x-82x^2+x^3} \quad ,
$$
then
$$
a_n^3 + b_n^3 = c_n^3 +(-1)^n \quad .
\eqno(1)
$$

This theorem is very attractive since it furnishes infinitely many {\it almost-counterexamples} to Fermat's last theorem for $n=3$.
It turns out that we have seen this identity before, so let's backtrack 25 years.

{\bf 1995:} As with many of Ramanujan's theorems, he did not give
 a proof, and left it to posterity (most notably Bruce Berndt and George Andrews and their students) to furnish proofs.
In 1995, Michael Hirschhorn gave a proof [Ha1] (see also [Ha2] and [HaHi]) and more interestingly, indicated
how Ramanujan may have discovered it. Hirschhorn's {\it exegesis} started with a {\bf polynomial identity}
that appeared in Ramanujan's notebooks.
$$
(A^2+7AB-9B^2)^3+ (2A^2-4AB+12B^2)^3 + (-2A^2-10B^2)^3+ (-A^2+9AB+B^2)^3 \, = \, 0 \quad.
\eqno(2)
$$
Hirschhorn speculated that this polynomial identity was the starting point of Ramanujan's discovery. Then he solved the 
quadratic diophantine equation
$$
A^2-9AB-B^2= \pm 1 \quad,
$$
in terms of the recursive sequence $h_n$ defined by $h_0=0,h_1=1$ and $h_{n+2}=9h_{n+1}+h_n$ and showed that
indeed it is satisfied by $A=h_{n+1}$ and $B=h_n$. Then Hirschhorn (and presumably Ramanujan) found generating functions
for the remaining quantities $A^2+7AB-9B^2$, $2A^2-4AB+12B^2$ and $2A^2+10B^2$. 

As noticed by Hirschhorn, once Ramanujan's amazing identity is conjectured, its formal proof is {\bf purely routine}.
In fact (see [Ha2]) checking it for the seven special cases $n=0,1,2,3,4,5,6$ consists a {\bf fully rigorous proof}.
As we will explain later, the {\bf broader context} is the $C$-finite ansatz ( Ch. 4 of [KP]; [Z2]).

But how in the world did Ramanujan come up with the equally amazing {\it polynomial} identity $(2)$? In this
case it is even more obvious that {\bf once conjectured}, the proof is a routine high-school-algebra computation, but 
how did Ramanujan find these four quadratic polynomials \hfill\break
$P_1(A,B),P_2(A,B),P_3(A,B),P_4(A,B)$ with integer coefficients such that
$$
P_1(A,B)^3+P_2(A,B)^3+P_3(A,B)^3+P_4(A,B)^3 \, = \, 0 \quad ?
$$
Now we have to backtrack further in time, about $55$ years ago.

{\bf 1965:} One of us (DZ) remembered that when he was fifteen-years-old  he was fascinated by  an article [J],
published in his then (and still now) favorite  journal ``gilyonot letmatemtika"  [in Hebrew]
by Eri Jabotinsky(1910-1969) in which he attempted to ``read Ramanujan's mind" and explain how he may have come up with such
infinite families of `sum of cubes identities' generalizing $3^3+4^3+5^3=6^3$.
The details were forgotten, but thanks to the efforts of Gadi Aleksandrowicz, from the Technion, Israel,
it has been archived and  made available here

{\tt https://gadial.github.io/netgar/} \quad .

We easily were  able to locate Jabotinsky's article, that so fascinated us $55$ years ago, and use it for this project. 
Jabotinsky's `trick' will be explained later.

{\bf Why this article}

Now that we know Ramanujan's `tricks' (at least as conjectured by Hirschhorn and Jabotinksy, and it does not matter whether
that's how Ramanujan actually did it), can they be used for other problems? Can they be turned into an algorithm that
can be taught to a computer? Using computer algebra, we were able to combine these two `tricks' 
and developed an algorithm to solve
diophantine equations of the more general form
$$
a \,X^3 + a \, Y^3 + b\, Z^3 \, = \, c \quad ,
$$
for {\it all} integers $a,b$ and a select set of $c$'s. 

{\bf Sample Theorem}: The front of this article  \hfill\break
{\tt https://sites.math.rutgers.edu/\~{}zeilberg/mamarim/mamarimhtml/ramacubic.html} \hfill\break
contains lots of theorems. For example the output file \hfill\break
{\tt https://sites.math.rutgers.edu/\~{}zeilberg/tokhniot/oRamanujanCubic1.txt} \hfill\break
contains $641$ theorems similar to Ramanujan's
theorem above. Let's just reproduce one here, as an appetizer.

{\bf Theorem (S.B. Ekhad)}: Define sequences of integers $a_n,b_n,c_n$, $n \geq 0$ , in terms of the following
generating functions
$$
\sum_{n \geq 0} a_n t^n \, = \,  {\frac {293155\,{t}^{2}+888826\,t-29}{ \left( 1-t \right)  \left( {t}^{2}-103682\,t+1 \right) }} \quad ,
$$
$$
\sum_{n \geq 0} b_n t^n \, = \, -{\frac {237169\,{t}^{2}+550798\,t+1}{ \left( 1-t \right)  \left( {t}^{2}-103682\,t+1 \right) }}
\quad ,
$$
$$
\sum_{n \geq 0} c_n t^n \, = \, {\frac {90601\,{t}^{2}-878594\,t+25}{ \left( 1-t \right)  \left( {t}^{2}-103682\,t+1 \right) }}
\quad ,
$$
then
$$
a_n^3 + 2 b_n^3 + 2 c_n^3 \, = \, 6859 \quad .
\eqno(1)
$$

{\bf The Maple packages} 

This article  is accompanied by three Maple packages.

$\bullet$  {\tt RamanujanCubic.txt} : This is the  main Maple package. It can automatically solve 
diophantine equation of the from $a\,X^3\,+ \, a\,Y^3 \,+\, b\,Z^3 \,=c\, $ for any given choice of integers $a,b$ and some $c$
(that depend on $a$ and $b$) by furnishing a {\it parametric solution} thereby furnishing {\it infinitely} many solutions.

$\bullet$ {\tt Pell.txt}: It solves  Pell, and Pell-like diophantine equations for quadratic binary forms,
using symbolic computation.

$\bullet$ {\tt RatDio.txt}: It uses the $C$-finite ansatz to manufacture complicated diophantine equations
that are known to have solutions beforehand, and also solves them in the quadratic case, using  the $C$-finite ansatz.

They are all available  from the front mentioned above, along with numerous sample input and output files.

{\bf Diophantine Equations In General}

Hilbert's dream, expressed in his 10th problem, to  find an algorithm that would decide whether  or not {\it any} diophantine equation is solvable
was famously shattered by Yuri Matiyashevich, standing on the shoulders of Julia Robinson, Martin Davis, and Hilary Putnam.
But for many {\it families}, e.g. Pell's equation $x^2\,-\,N\,y^2=1$ (for any non-square) $N$ , and Fermat's 
$x^n+y^n-z^n=0$ ($n>2$) it was possible to {\it decide}, positively and negatively, respectively. But let us remark that by
going {\it backwards}, and with the aid of computer algebra, we can easily {\bf concoct} many diophantine equations
that are {\it guaranteed} to have infinitely many solutions.

Recall that, famously, the equation $x^2+y^2=z^2$ has the {\it parametric} solution
$$
x=m^2-n^2 \quad, \quad y=2mn \quad, \quad z=m^2+n^2 \quad .
$$

Let's start with any polynomials (with integer coefficients) $P(m,n)$, $Q(m,n)$ and $R(m,n)$,  and define
$$
x=P(m,n) \quad, \quad y=Q(m,n) \quad, \quad z=R(m,n) \quad .
$$
{\it Eliminating} $m,n$ from these equations (using, most efficiently, the Buchberger algorithm) we immediately find
a polynomial (with integer coefficients) $S(x,y,z)$ such that $S(P(m,n),Q(m,n),R(m,n))$ is identically $0$. So we can cook-up many theorems.

Suppose that, in  a {\it counterfactual world},
you never heard of Pythagorean triples, 
and defined, out of the blue, $x=m^2-n^2 \quad, \quad y=2mn \quad, \quad z=m^2+n^2 \quad$, and  wondered how $x,y,z$ are related
{\it implicitly}. If you enter, in Maple:

{\tt Groebner[Basis]({x-(m**2-n**2),y-2*m*n,z-m**2-n**2},plex(m,n,x,y,z))[1];} \quad ,

you would get the  beautiful relationship
$$
{x}^{2}+{y}^{2}-{z}^{2}=0 \quad .
$$

Typing instead

{\tt Groebner[Basis]({x-(2*m**2-3*n**2),y-2*m*n,z-m**2-n**2},plex(m,n,x,y,z))[1];}

gives the  less good-looking relationship
$$
4\,{x}^{2}+4\,xz+25\,{y}^{2}-24\,{z}^{2} \, = \, 0 \quad,
$$

but you can publish a paper with the next deep theorem.

{\bf Theorem}: The diophantine equation
$$
4\,{x}^{2}+4\,xz+25\,{y}^{2}-24\,{z}^{2} \, = \, 0  \quad ,
$$
has (doubly) infinite many solutions given by the parametric equation
$$
x=2m^2-3n^2 \quad , \quad  y=2mn \quad, \quad z=m^2+n^2 \quad  .
$$
Proof: (check!).

This way one can make up many such deep theorems. Here is another example

{\tt Groebner[Basis]({x-(m**3-n**3),y-m**2*n-m*n**2,z-(m**3+n**3)},plex(m,n,x,y,z))[1];} gives

$$
3\,{x}^{2}y+{x}^{2}z+4\,{y}^{3}-3\,y{z}^{2}-{z}^{3} \quad,
$$

and we immediately have the following deep theorem.

{\bf Theorem}: The diophantine equation
$$
3\,{x}^{2}y+{x}^{2}z+4\,{y}^{3}-3\,y{z}^{2}-{z}^{3} \,  = 0 \, \quad,
$$
has doubly-infinite many solutions given by the parametric equation
$$
x=m^3-n^3 \quad , \quad  y=m^2n+mn^2 \quad, \quad z=m^3+n^3 \quad   .
$$

Proof: (check!).

You are welcome to create your own deep theorems, for which the equation has (doubly!) infinitely many solutions.

To create random examples for which we know a-priori that there are {\bf no} solutions is  a bit harder, but you can "cheat" and
take any diophantine equation, e.g. $x^3+y^3+z^3=0$ for which it has been already proven that there exist no (non-trivial) solutions.
Then replace $x$, $y$ and $z$ by arbitrary expressions. For example, typing, in a Maple session,

{\tt expand(subs({x = 6*x+7*y-9*z, y = 6*x-5*y+4*z, z = -8*x-3*y+3*z},x**3+y**3+z**3));}

gives
$$
-80\,{x}^{3}-360\,{x}^{2}y+36\,{x}^{2}z+1116\,x{y}^{2}-2556\,xyz+1530\,x{z}^{2}+191\,{y}^{3}-942\,{y}^{2}z+1380\,y{z}^{2}-638\,{z}^{3}
$$

and we immediately have

{\bf Theorem}: The diophantine equation
$$
-80\,{x}^{3}-360\,{x}^{2}y+36\,{x}^{2}z+1116\,x{y}^{2}-2556\,xyz+1530\,x{z}^{2}+191\,{y}^{3}-942\,{y}^{2}z+1380\,y{z}^{2}-638\,{z}^{3} \, = \, 0
$$
has no (non-trivial) solutions.

{\bf Creating Diophantine Equations using the $C$-finite ansatz}

Another way to  {\bf start with the answer} is to use the $C$-finite ansatz ([KP],ch. 4; [Z2]).

Let's start with quadratic equations. We have the following theorem, that once stated is easy to prove in the
context of the $C$-finite ansatz.

{\bf General Theorem:} Let $c0,c1$, $d0,d1$ and $k$ be integers and define two sequences of integers $a(n),b(n)$
$$
\sum_{n=0}^{\infty} \, a(n) x^n \, = \, \frac{c0+c1x}{1-kx+x^2} \quad,
$$
$$
\sum_{n=0}^{\infty} \, b(n) x^n \, = \, \frac{d0+d1x}{1-kx+x^2} \quad,
$$
then there exists a homogeneous binary quadratic form $P(X,Y)$, with integer coefficients, and an integer $C$, such that for all $n \geq 0$
$$
P( a(n), b(n) ) \, = \, C \quad .       
$$

In fact, more explicitly
$$
 \left( {\it c0}\,{\it d1}-{\it c1}\,{\it d0} \right)^{2} \left( {\it d0}\,{\it d1}\,k+{{\it d0}}^{2}+{{\it d1}}^{2} \right) a(n)^2
$$
$$
- \left( {\it c0}\,{\it d1}-{\it c1}\,{\it d0} \right)^{2} \left( {\it c0}\,{\it d1}\,k+{\it c1}\,{\it d0}\,k+2\,{\it c0}\,{\it d0}+2\,{\it c1}\,{\it d1}
 \right) a(n)\,b(n)
$$
$$
+
 \left( {\it c0}\,{\it d1}-{\it c1}\,{\it d0} \right)^{2} \left( {\it c0}\,{\it c1}\,k+{{\it c0}}^{2}+{{\it c1}}^{2} \right)\,b(n)^2
\, = \,
 \left( {\it c0}\,{\it d1}-{\it c1}\,{\it d0} \right)^{4} \quad .
$$

An important special case solves {\bf  Pell's Equation}.

{\bf Special Theorem (solving Pell's Equation)} Let $b$ and $k$ be integers. Define two sequences of integers $A(n),B(n)$
$$
\sum_{n=0}^{\infty} \, A(n) x^n \, = \, \frac{1-kx}{1-2kx+x^2} \quad,
$$
$$
\sum_{n=0}^{\infty} \, B(n) x^n \, = \, \frac{bx}{1-2kx+x^2} \quad .
$$
Let
$$
P(X,Y)=X^2- \left ( \frac{k^2-1}{b^2} \right ) Y^2-1 \quad , 
$$
then
$$
P( A(n), B(n) ) \, = \, 0 \quad .       
$$
In particular, once we found the smallest non-trivial solution, $(k,b)$,  of {\bf Pell's equation}, $k^2\,-\, N\,b^2 \, =\,1$, 
it enables us to get
{\bf all} solutions. Of course this is equivalent to the standard way, using  expressions of the form
$(k+\sqrt{N} b)^n$, and extracting the coefficients of $\sqrt{N}$, and $1$, but our approach is simpler, and puts it in the
context of solving diophantine equations within the C-finite ansatz.

We also have

{\bf Another General Theorem:} Let $c_0,c_1$, $d_0,d_1$ and $k$ be integers and define two sequences of integers $a(n),b(n)$
$$
\sum_{n=0}^{\infty} \, a(n) x^n \, = \, \frac{c_0+c_1x}{1-kx-x^2} \quad,
$$
$$
\sum_{n=0}^{\infty} \, b(n) x^n \, = \, \frac{d_0+d_1x}{1-kx-x^2} \quad,
$$
then there exists a homogeneous binary quadratic from $P(X,Y)$ with integer coefficients, and an integer $C$ such that for all $n \geq 0$
$$
P( a(n), b(n) ) \, = \, C(-1)^n \quad .       
$$

Even more generally, we have the following theorem. 

{\bf Very General Theorem:} Let $d$ be an integer larger than $1$, and let $c_{i,j}$ $1 \leq i \leq d, 0 \leq j \leq d-1$ and
$k_1 , \dots,  k_{d-1}$ be integers. Define $d$ integer sequences $a_i(n)$ ($1 \leq i \leq d$) by the generating functions
$$
\sum_{n=0}^{\infty} \, a_i(n) x^n \, = \, \frac{c_{i,0} + \cdots + c_{i,d-1} x^{d-1} }{1-k_1x-k_2 x^2 - \dots k_{d-1} x^{d-1}+(-1)^d x^d} \quad, \quad
(1 \leq i \leq d) \quad ,
$$
then there exists a homogeneous polynomial $P(X_1, \dots, X_d)$, with integer coefficients, of degree $d$ and an integer $C$ such that
$$
P( a_1(n), a_2(n), \dots, a_d(n) ) \, = \, C \quad .       
$$

Analogously

{\bf Very General Theorem':} Let $d$ be an integer larger than $1$, and let $c_{i,j}$ $1 \leq i \leq d, 0 \leq j \leq d-1$ and
$k_1 , \dots,  k_{d-1}$ be integers. Define $d$ integer sequences $a_i(n)$ ($1 \leq i \leq d$) by the generating functions
$$
\sum_{n=0}^{\infty} \, a_i(n) x^n \, = \, \frac{c_{i,0} + \cdots + c_{i,d-1} x^{d-1} }{1-k_1x-k_2 x^2 - \dots k_{d-1} x^{d-1}+(-1)^{d+1} x^d} \quad,
(1 \leq i \leq d)
$$
then there exists a homogeneous polynomial $P(X_1, \dots, X_d)$ of degree $d$ (with integer coefficients) and an integer $C$ such that
$$
P( a_1(n), a_2(n), \dots, a_d(n) ) \, = \, C (-1)^n\quad .       
$$

{\bf Back to the Quadratic Case}

Using symbolic computations and experimental mathematics, we can get an alternative approach for {\it solving}
quadratic diophantine equations. Of course, this is all classical, and  modern treatments can be found, e.g.
in the excellent texts [S] and [AA]. Our approach also uses continued fractions, but in a much more
simple-minded way, without the usual human-generated infra-structure. This is implemented in the Maple package
{\tt Pell.txt}. The main procedure is  {\tt SolQuad(Q,m,n,t,K)} that inputs a quadratic binary form with integer coefficients, $Q(m,n)$, in the
variables $m$ and $n$ and a symbol $t$ and a positive integer $K$ (a parameter for `guessing' it can always be increased
until success is reached), and outputs two rational functions, $f_1(t),f_2(t)$  (with identical, quadratic, monic denominators) that are the generating functions
for the sequences $a(i),b(i)$ that satisfy $Q(a(i),b(i))=constant$. From now on we can take it as a {\bf black box}.

{\bf Eri Jabotinsky's trick}

Ramanujan was not the first one to have a parametric solution to the equation $X^3+Y^3+Z^3+W^3=0$, this
honor goes to Euler (see [D], vol. 2, Chap. XXI, pp. 552-553), but Ramanujan had quite a few of them,
and probably found them all by himself. In a delightful article meant for teenagers, Eri Jabotinsky [J] explained
how Ramanujan {\it might} have found them.

Suppose that you have two distinct solutions $(x,y,z,w)$ and $(x',y',z',w')$ of the equation
$X^3+Y^3+Z^3+W^3=0$, in other words
$$
x^3+y^3+z^3+w^3=0 \quad, \quad (x')^3+(y')^3+(z')^3+ (w')^3=0 \quad .
$$

Let $c$ and $d$ be {\it undetermined coefficients} for now, and let's try to find a brand new solution  of the form
$$
X=cx+dx' \quad, \quad 
Y=cy+dy' \quad, \quad 
Z=cz+dz' \quad, \quad 
W=cw+dw' \quad .
$$
Expanding $X^3+Y^3+Z^3+W^3$ we get
$$
c^3 x^3 + 3 c^2 d x^2 x' + 3 cd^2 x x'^2 + d^3 x'^3 \quad +
$$
$$
c^3 y^3 + 3 c^2 d y^2 y' + 3 cd^2 y y'^2 + d^3 y'^3 \quad +
$$
$$
c^3 z^3 + 3 c^2 d z^2 z' + 3 cd^2 z z'^2 + d^3 z'^3 \quad +
$$
$$
c^3 w^3 + 3 c^2 d w^2 w' + 3 cd^2 w w'^2 + d^3 w'^3 \quad .
$$
Collecting terms, and equating to $0$, we get 
$$
c^3 (x^3+y^3+z^3+w^3) 
$$
$$
+ 3cd  ( c (x^2 x'+ y^2 y' + z^2 z' + w^2 w') + d (x x'^2+ y y'^2 + z z'^2 + w w'^2) )
$$
$$
+d^3 (x'^3+y'^3+z'^3+w'^3)  \, = \, 0 \quad .
$$
The first and last terms vanish by assumption, so we need to choose integers $c$ and $d$ such that
$$
 c (x^2 x'+ y^2 y' + z^2 z' + w^2 w') + d (x x'^2+ y y'^2 + z z'^2 + w w'^2) \, =0 .
$$
Taking
$$
c=x'^2 x+ y'^2 y + z'^2 z + w'^2 w \quad, \quad
d=-(x^2 x'+ y^2 y' + z^2 z' + w^2 w') \quad ,
$$
would work, so out of the original two {\bf  numerical} special solutions $(x,y,z,w)$ and $(x',y',z',w')$ we get yet-another
specific numerical solution. By pure inspection, (even by hand!) it is very fast to come up with specific
numerical solutions, the most famous ones being $(3,4,5,-6)$ and $(9,10,-1,-12)$ [of taxicab fame].
But in addition, there is a trivial `doubly-infinite'
solution $(m,-m,n,-n)$. By using the Eri Jabotinsky trick (that Jabotinsky conjectured was Ramanujan's way),
`morphing' any specific solution with the symbolic solution $(m,-m,n,-n)$, we get a doubly-infinite {\bf non-trivial}
solution consisting of quadratic polynomials in $m$ and $n$.

{\bf  Back to Mike Hirschhorn}

Mike Hirshhorn's starting point (and presumably Ramanujan's) was  the identity
$$
(A^2+7AB-9B^2)^3+ (2A^2-4AB+12B^2)^3 + (-2A^2-10B^2)^3+ (-A^2+9AB+B^2)^3 \, = \, 0 \quad.
\eqno(2)
$$

Using human ingenuity, he then solved the diophantine equation  $-A^2+9AB+B^2 \, = \, \pm 1$, by [essentially]
finding two rational functions whose respective Taylor coefficients satisfy this equation. Then
using further {\it manipulatorics}, he  found generating functions for the other quantities,
thereby rediscovering, {\it ab initio},  Ramanujan's theorem. 
All this can be streamlined in the context of the $C$-finite  ansatz, and has been implemented in the Maple package
{\tt RamanujanCubic.txt}, mentioned above. 
But what is the point? Who cares about a computerized redux of the beautiful human-generated work of Hirschhorn,
Jabotinsky  (and Ramanujan).

The point is that we can {\bf generalize}. Eri Jabotinsky's trick works just as well for
the more general cubic diophantine equation
$$
a\,X^3 \, + \, a \,Y^3 \,+\,  b \,Z^3 \,+\,  b\,W^3 \, = \, 0,
$$
for {\it any} integers $a$ and $b$.
It is easy, by pure brute force, to come up with specific numerical solutions (even by hand). The symbolic solution 
$(m,-m,n,-n)$ is still a solution. Using the Jabotinsky trick, we get
 four quadratic polynomials $P_1(m,n),P_2(m,n),P_3(m,n),P_4(m,n)$ such that
$$
a\,P_1(m,n)^3 \, +\, a\,P_2(m,n)^3 \, + \, b\,P_3(m,n)^3 \,+ \, b P_4(m,n)^3 \, = \, 0 \quad .
$$
Then we can use (our version of) the well-known algorithm (implemented in our Maple package) to find a pair of rational  functions $f_1(t),f_2(t)$ 
whose respective Taylor coefficients solve $P_4(m,n)=1$ (or failing this [there is not always a solution])
$P_4(m,n)=c$ for some fixed small integer $c$. Then using the `$C$-finite  calculator' [Z2], also built-in into the package, one gets
as many Ramanujan-like theorems solving diophantine equations of the form $a\, X^3\,+\,a \,Y^3 \, +b\, \,Z^3=c$, (with a different $c$, of course)
for any desired $a$ and $b$, and for  some emerging constant $c$, like the sample theorem above.

For $640$ additional such theorems, see the output file

{\tt https://sites.math.rutgers.edu/\~{}zeilberg/tokhniot/oRamanujanCubic1.txt} \quad .

Of course, using the Maple package 

{\tt http://www.math.rutgers.edu/\~{}zeilberg/tokhniot/RamanujanCubic.txt} \quad ,

readers can generate many more such theorems.

{\bf Morals}

{\bf 1.} Read Children books, especially such delightful ones as Amy Alznauer's about Ramanujan.

{\bf 2.} Reread your favorite math articles from when you were a teenager, in this case made possible by
 Gadi Aleksandrowicz. 

{\bf 3.} Reread insightful papers, like [Ha1] trying to do a {\bf human deconstruction} of geniuses like Ramanujan.

{\bf 4.} Teach these tricks to a computer, doing  {\bf computerized deconstruction} (see e.g. [Z1])

{\bf 5.} By doing (usually) minor tweaking to the program, generalize it, and generate  as many Ramanujan-like
theorems as you wish.

{\bf References}

[AA] Titu Andreescu and Dorin Andrica, {``Quadratic Diophantine equations''},
 Developments in Mathematics {\bf 40}. Springer, New York, 2015.

[Al] Amy Alznauer (author) and Daniel Miyares (illustrator), {\it ``The  Boy Who Dreamed of Infinity''},
Candlewick Press, 2020

[D] Leonard E. Dickson, {``History of the Theory of Numbers''} ($2$ volumes),
Carnegie Institute, 1919, 1920, and 1923. Reprinted by Chelsea Publishing Company, 1992.

[HaHi] Jung Hun Han and Michael D. Hirschhorn, {\it Another look at an amazing identity of Ramanujan},
Mathematics Magazine {\bf 79} (2006), 302-304.

[Ha1] Michael D. Hirschhorn, {\it An amazing identity of Ramanujan}, Mathematics Magazine {\bf 68} (1995), 199-201.

[Ha2] Michael D. Hirschhorn, {\it A proof in the spirit of Zeilberger of an amazing identity of Ramanujan}, 
Mathematics Magazine {\bf 69} (1996), 267-269.

\vfill\break
[J] Eri Jabotinsky,  {\it The Pythagorean numbers and the cubic numbers that resemble them} [in Hebrew],
{\it dapim le-matematika ulefisika} [Pages for mathematics and physics] vol. 2, issue 1 [Feb. 1946].
[edited by Dov  Jarden and Dr. Theodore Motzkin\footnote{${}^1$}
{\eightrm  Eri Jabotinsky and Theodore Motzkin were both sons of Zionist leaders, Ze'ev Jabotinsky
[the leader of the militant branch that inspired the Jewish terrorist groups Etzel and Lechi], and
Leo Motzkin, a close friend of Theodore Herzl. Other examples are Herb Wilf, whose father Alexander 
lead the fund-raising drive for the {\it irgun} before 1948, and Gil Kalai, whose father, Hanoch, was second-in-command in Yair Stern's
organization, {\it Lechi}.
}, 
published by  Eri Jabotinsky (publisher) Ltd.] 
Downloadable from \hfill\break
{\tt https://gadial.github.io/netgar/} \quad .
%{\tt https://docs.google.com/file/d/0B-\_8w6IKpNuUc3ZzMXk5ZnNDbFE/view} \quad . \hfill\break
Reprinted in {\it gilyonot lematematika} (with a forward by Joseph Gillis) ca. 1965.

[KP] Manuel Kauers  and Peter Paule, {\it ``The Concrete Tetrahedron''}, Springer, 2011.

[R] Srinivasa Ramanujan,{\it ``The Lost Notebook and Other Unpublished papers''}, Narosa, New Delhi, 1988.

[S] Nigel P. Smart, {\it `` The Algorithmic Resolution of  Diophantine Equations''}, London Mathematical Society Student Texts {\bf 41},
Cambridge University Press, 1998.

[Z1] Doron Zeilberger, {\it Computerized Deconstruction},  Adv. Appl. Math. {\bf 30} (2003), 633-654. \hfill\break
{\tt https://sites.math.rutgers.edu/\~{}zeilberg/mamarim/mamarimhtml/derrida.html} \quad .

[Z2] Doron Zeilberger, {\it The C-finite Ansatz},  Ramanujan Journal {\bf 31} (2013), 23-32. \hfill\break
{\tt https://sites.math.rutgers.edu/\~{}zeilberg/mamarim/mamarimhtml/cfinite.html} \quad .

\bigskip
\bigskip
\hrule
\bigskip
Shalosh B. Ekhad, c/o D. Zeilberger, Department of Mathematics, Rutgers University (New Brunswick), Hill Center-Busch Campus, 110 Frelinghuysen
Rd., Piscataway, NJ 08854-8019, USA. \hfill\break
Email: {\tt ShaloshBEkhad at gmail dot com}   \quad .
\bigskip
Doron Zeilberger, Department of Mathematics, Rutgers University (New Brunswick), Hill Center-Busch Campus, 110 Frelinghuysen
Rd., Piscataway, NJ 08854-8019, USA. \hfill\break
Email: {\tt DoronZeil at gmail  dot com}   \quad .
\bigskip
\hrule
\bigskip
Written: July 30, 2020.
\end